\documentclass[11pt]{article}
\renewcommand{\baselinestretch}{1.2}
\newcommand{\dated}{\mbox{} \hfill {\small [{\tt \today}]}} \usepackage{amsmath,amssymb,amsfonts}
\newcommand{\comps}{{\mathbb C}}

\newcommand{\ints}{{\mathbb Z}}

\newcommand{\torus}{{\mathbb T}}

\newcommand{\cstar}{{C^\ast}}

\newcommand{\A}{{\mathfrak A}}

\usepackage{amsthm,enumerate}
\newtheorem*{theorem}{Theorem}

\title{Why Banach algebras?}
\author{\textit{Volker Runde}}
\date{}
\begin{document}
\maketitle
Most readers probably know what a Banach algebra is, but to bring everyone to the same page, here is the definition: A \emph{Banach algebra} is an algebra---always meaning: linear, associative, and not necessarily unital---$\A$ over $\comps$ which is also a Banach space such that
\begin{equation} \label{submult} \tag{\mbox{$\ast$}}
   \| ab \| \leq \| a \| \| b \| \qquad (a \in \A)
\end{equation}
holds. 
\par
The question in the title is deliberately vague. It is, in fact, shorthand for several questions.
\par
\subsubsection*{Why are they called Banach algebras?} Like the Banach spaces, they are named after the Polish mathematician Stefan Banach (1892--1945) who had introduced the concept of Banach space, but Banach had never studied Banach algebras: they are named after him simply because they are algebras that happen to be Banach spaces. The first to have explicitly defined them---under the name ``linear metric rings''---seems to have been the Japanese mathematician Mitio Nagumo in 1936 (\cite{Nagumo}), but the mathematician whose work really got the theory of Banach algebras on its way was Israel M.\ Gelfand (1913--2009), who introduced them under the name ``normed rings'' (\cite{Gelfand}). In the classical monograph \cite{BonsallDuncan}, the authors write that, if they had it their way, they would rather speak of ``Gelfand algebras''. I agree. But in 1945 , Warren Ambrose (1914--1995) came up with the name ``Banach algebras'' (\cite{Ambrose}) , and it has stuck ever since.
\par 
\subsubsection*{Why am I working on Banach algebras?} It just happened. I started studying mathematics in 1984 at the University of M\"unster in Germany, and American expatriate George Maltese (1931--2009) was teaching the introductory analysis sequence. As he was an excellent teacher, I also took his subsequent courses on functional analysis and, eventually, on Banach algebras. I wrote my ``Diplomarbeit'' under his supervision on a topic from the theory of Banach algebras, and the area has kept me hooked ever since.
\par 
This, of course, begets the question of what has been keeping me hooked on Banach algebras for now over two decades. The answer is multi-layered and necessarily dependent on my personal mathematical tastes.
\par 
First of all, Banach algebra theory is elegant (at least most of it). 
\par
Commutative Banach algebras allow for a powerful representation theory, which is now known as \emph{Gelfand theory} after its creator. Suppose that $\A$ is a commutative Banach algebra (with an identity element $e$, for the sake of simplicity), and let
\[
  \Phi_\A := \{ \phi \!: \A \to \comps : \text{$\phi$ is non-zero, linear, and multiplicative} \}.
\]
It is easy to see that $\Phi_\A$ lies is the unit ball of the dual space $\A^\ast$, and restricting the weak$^\ast$ topology of $\A^\ast$ to $\Phi_\A$ turns it into a compact Hausdorff space. For each $a \in \A$, we thus obtain a continuous function $\hat{a}$ on $\Phi_\A$ by letting $\hat{a}(\phi) := \langle a, \phi \rangle$ for $\phi \in \Phi_\A$; the function $\hat{a}$ is called the \emph{Gelfand transform} of $a$. The remarkable insight of Gelfand was that $a$ is invertible in $\A$ if and only if $\hat{a}$ is invertible as a continuous function on $\Phi_\A$, i.e., has no zeros. 
\par 
The Gelfand transform can be thought of as an abstract Fourier transform, and it was to Fourier analysis where the still young theory of Banach algebra had its first striking and strikingly beautiful application. In \cite{Wiener}, Norbert Wiener (1894--1964) proved:
\begin{theorem}
Let $\torus$ be the unit circle in $\comps$ and let $f \!: \torus \to \comps$ be a continuous function without zeros whose Fourier series $\sum_{n \in \ints} \hat{f}(n)$ converges absolutely. Then the Fourier series of $1/f$ converges absolutely as well.
\end{theorem}
\par 
Wiener's proof in \cite{Wiener} was hard Fourier analysis. Using his representation theory, Gelfand gave an astonishingly short and simple proof. Let
\[
  \mathcal{W} := \{ f \!: \torus \to \comps : \text{$f$ is continuous with absolutely converging Fourier series} \},
\]
and set $\| f \| := \sum_{n \in \ints} | \hat{f}(n)|$ for $f \in \mathcal{W}$. It is easy to see that this turns $\mathcal{W}$ into a commutative Banach algebra, and it is equally simple to see that $\Phi_\A$ can be identified with $\torus$ in such a way that $\hat{f} = f$. Thus, if $f \in \mathcal{W}$ has no zeros, it must already be invertible in $\mathcal{W}$, i.e., $1/f \in \mathcal{W}$. 
\par
Secondly, there is more to the axioms of a Banach algebra than meets the eye at first glance. We have an algebraic structure---an algebra---and an analytic structure---a Banach space---, and the two structures are linked by means of the inequality (\ref{submult}), which guarantees that multiplication in Banach algebras is continuous. As it turns out, the relationship between the analytic and the algebraic structure is much more subtle.
\par 
If we have an algebra equipped with a norm $\| \cdot \|$ turning it into a Banach algebra are there other norms---meaning: not equivalent to $\| \cdot \|$---that accomplish the same? In general, the answer is no: take a linear space which is a Banach space under two inequivalent norms and turn it into a Banach algebra by defining the product of any two elements to be zero. But there are situations where a given algebra admits---up to equivalence---only one norm turning it into a Banach algebra. For instance, let $K$ be a compact Hausdorff space. Then $\A = \mathcal{C}(K)$, the algebra of all continuous functions on $K$, is a Banach algebra under the supremum norm $\| \cdot \|_\infty$. A simple closed graph argument shows that any norm turning $\A$ into a Banach algebras is already equivalent to $\| \cdot \|_\infty$. This is only a very easy special case of the following famous theorem (\cite{Johnson}) due to Barry E.\ Johnson (1937--2002):
\begin{theorem}[Johnson's uniqueness of norm theorem]
Let $\A$ be a semisimple---meaning: with zero Jacobson radical---Banach algebra. Then any norm on $\A$ turning it into a Banach algebra is equivalent to the given norm. 
\end{theorem}
\par 
An algebraic hypothesis like semisimplicity thus forces the analytic structure to be unique.
\par 
In \cite{Kaplansky}, Irving Kaplansky (1917--2006) asked if, given a compact Hausdorff space $K$, any submultiplicative, i.e., satisfying (\ref{submult}), but not a priori complete norm on $\A = \mathcal{C}(K)$, was already equivalent to $\| \cdot \|_\infty$. This doesn't look like a hard problem, but it was, and it took almost three decades to solve it (\cite{DalesEsterle}):
\begin{theorem}
Let $K$ be an infinite compact Hausdorff, \emph{and suppose that the continuum hypothesis holds}. Then there is a submultiplicative norm on $\mathcal{C}(K)$ not equivalent to $\| \cdot \|_\infty$.
\end{theorem}
\par 
Strangely, the demand that the continuum hypothesis hold cannot be dropped (\cite{DalesWoodin}).
\subsubsection*{Why should the rest of mathematics care about Banach algebras?}
Banach algebras show up naturally in many areas of analysis: If $E$ is a Banach space, then the bounded linear operators on $E$ form a Banach algebra, as do the compact operators, or the nuclear ones, etc.; if $G$ is a locally compact group, then the space $L^1(G)$ of integrable functions (with respect to left Haar measure) is a Banach algebra under under convolution; and finally, even though the area of operator algebras has developed a lift pretty much independent of general Banach algebra theory, all $\cstar$- and von Neumann algebras are, first of all, Banach algebras. If we have an analytic object that has a Banach algebra naturally associated with it, then this algebra can provide us with further insight into the nature of the underlying object.
\par 
Recall that a locally compact group $G$ is called \emph{amenable} if it has a \emph{left invariant mean}, i.e., a bounded linear functional $M \!: L^\infty(G) \to \comps$ with $\| M \| = 1 = \langle 1, M \rangle$ such that
\[
  \langle L_x \phi, M \rangle = \langle \phi, M \rangle \qquad (\phi \in L^\infty(G), \, x \in G);
\]
here, $L_x \phi$ stands for the left translate of $\phi$ by $x$. Abelian groups are amenable, and so are compact groups, but the free group on two generators isn't; also, amenability is stable under standard constructions such as subgroups, quotients, and extensions. The concept goes back to John von Neumann (1903--1957) (\cite{vonNeumann}).
\par
There is also the notion of an amenable Banach algebra; it is due to Johnson (\cite{Johnsonmemoir}). There are various equivalent characterizations of amenable Banach algebras (see \cite{Runde}). The one that probably requires the least background is as follows. Given a Banach algebra $\A$, a bimodule $E$ over $\A$ which is also a Banach space such that the module actions of $\A$ on $E$ are continuous is called a \emph{Banach $\A$-bimodule}. We call $\A$ \emph{amenable} if, for every Banach $\A$-bimodule $E$ and for every derivation $D \!: \A \to E$, i.e., a bounded linear map satisfying $D(ab) = a \cdot Db + (Da) \cdot b$ for  $a, b \in \A$, there is a bounded net $( x_\alpha )_\alpha$ in $E$ such that
\[
  Da = \lim_\alpha a \cdot x_\alpha - x_\alpha \cdot a \qquad (a \in \A).
\]
\par 
The choice of the adjective ``amenable''---both for certain locally compact groups and for a class of Banach algebras---suggests that there is a link between the two, and indeed, there is one (\cite[Theorem 2.5]{Johnsonmemoir}):
\begin{theorem}
The following are equivalent for a locally compact group $G$:
\begin{enumerate}[\rm (i)]
\item $G$ is amenable;
\item $L^1(G)$ is an amenable Banach algebra.
\end{enumerate}
\end{theorem}
\par
This means that the amenability of a locally compact group can be capture through an entirely Banach algebraic property of its $L^1$-algebra.
\par 
This leads to an obvious question: What does amenability mean for other classes of Banach algebras, say for $\cstar$-algebras or algebras of operators on Banach spaces?
\par 
Already in \cite{Johnsonmemoir}, Johnson asked if every $\cstar$-algebra was amenable. The question soon turned out to have a negative answer, but that left the problem to characterize those $\cstar$-algebra that are amenable Banach algebras. It took the effort of a number of ma\-the\-ma\-ti\-cians---among them Alain Connes (\cite{Connes}) and Uffe Haagerup (\cite{Haagerup})---to obtain such a characterization:
\begin{theorem}
A $\cstar$-algebra is amenable if and only if it is nuclear.
\end{theorem}
\par 
Roughly speaking, a $\cstar$-algebra is nuclear if the identity map on it factors asymp\-to\-ti\-cal\-ly---in a sense that I won't make precise here---through full matrix algebras; for details, see \cite{Takesaki}.
\par
Another question from \cite{Johnsonmemoir} was if $\mathcal{B}(E)$, the Banach algebra of all bounded linear operators on a Banach space $E$, could be amenable for infinite-dimensional $E$. As amenable Banach algebras have a tendency to be ``small''---whatever that may mean precisely---, the knee jerk response would be a clear ``no''. Nevertheless, as a surprising by-product of the recent solution of the ``scalar plus compact problem'' (\cite{ArgyrosHaydon}), there are infinite-dimensional Banach spaces $E$ with $\mathcal{B}(E)$ amenable. Still, that example ought to be an exception rather than the rule, and indeed, only a few years after the publication of \cite{Johnsonmemoir}, it had become clear that $\mathcal{B}(\ell^2)$ is not amenable (\cite{Wassermann}). Of course, this strongly suggests that $\mathcal{B}(\ell^p)$ is not amenable for \emph{any} $p \in [1,\infty]$, but until the early 21$^\mathrm{st}$ there was no proof for the non-amenability of $\mathcal{B}(\ell^p)$ except for $p=2$. Then Charles J.\ Read proved that $\mathcal{B}(\ell^1)$ is not amenable (\cite{Read}). Simplifications of his proof (\cite{Pisier} and \cite{Ozawa}) led to a proof that $\mathcal{B}(\ell^\infty)$ is not amenable, and eventually the case for $p \in (1,\infty)$ was settled (\cite{RundeJAMS}). We thus have:
\begin{theorem}
Let $p \in [1,\infty]$. Then $\mathcal{B}(\ell^p)$ is not amenable.
\end{theorem}
\renewcommand{\baselinestretch}{1.0}
\dated
\vfill
\begin{tabbing}
\textit{Author's address}: \= Department of Mathematical and Statistical Sciences \\
\> University of Alberta \\
\> Edmonton, Alberta \\
\> Canada T6G 2G1 \\[\medskipamount]
\textit{E-mail}: \> \texttt{vrunde@ualberta.ca} \\[\medskipamount]
\textit{URL}: \> \texttt{http://www.math.ualberta.ca/$^\sim$runde/}
\end{tabbing}          
\end{document}